\documentclass[12pt]{article}
\usepackage{amsmath,latexsym,amssymb}
\usepackage{enumerate}
\usepackage{amsthm}
\usepackage{pgfplots}

\usepgfplotslibrary{fillbetween,decorations.softclip}
\pgfplotsset{compat=1.14}

\newtheorem{theorem}{Theorem}[section]

\newtheorem{lemma}[theorem]{Lemma}

\newtheorem{proposition}[theorem]{Proposition}
\newtheorem{conjecture}[theorem]{Conjecture}

\newcommand{\address}{Address: Facultad de Ciencias F\'isico-Matem\'aticas, Universidad Aut\'onoma de Sinaloa, Av. de las Am\'ericas y Blvd. Universitarios S/N Culiac\'an, Sinaloa, M\'exico C.P.80010; E-mail: islas.jose@uas.edu.mx}

\title{ Oscillation criteria for stopping near the top of a random walk}
\author{Jos\'e A. Islas\footnote{\address}}

\begin{document}

\maketitle

\begin{abstract}
Consider the problem of maximizing the probability of stopping with one of the two highest values in a Bernoulli random walk with arbitrary parameter $p$ and finite time horizon $n$.  Allaart \cite{Allaart} proved that the optimal strategy is determined by an interesting sequence of constants $\{p_{n}\}$. He conjectured the asymptotic behavior to be $1/2$. In this work the best lower bound for this sequence is found and more of its properties are proven towards solving the conjecture.

\end{abstract}

\section{Introduction}

Consider $\{S_n\}_{n=0,1...}$ to be a Bernoulli random walk with parameter $p \in (0,1)$. This is $S_0:=0$ and for $n\geq1$, $S_n:=X_1+...+X_n$, where $X_1, X_2,...$ are independent, identically distributed random variables with $P(X_1=1)=p$ and $P(X_1=-1)=q:=1-p$. Let $M_n:=\max\{S_0,S_1,...,S_n\}$, $n \in \mathbb{N}$. For a given finite time $N$, what is the maximum probability of ``stopping at the top'' of a random walk, that is, what is the stopping time $\tau$ adapted to the process $\{S_n\}$ that maximizes $P(S_{\tau}=M_N)$? 
This problem (also called the classical secretary problem with a random walk) was examined by Hlynka and Sheahan \cite{Hlynka} for $p=1/2$. Later Yam et al. \cite{Yam} considered the general case and found that the optimal strategy is simple: If $p>1/2$, the rule $\tau\equiv N$ is the unique optimal stopping rule; if $p<1/2$, $\tau\equiv 0$ is the unique optimal stopping rule; and if $p=1/2$, any rule $\tau$ such that $P(S_{\tau}=M_{\tau} $ or $ \tau=N)=1$ is optimal. This problem was extended by Allaart \cite{Allaart} to maximize the probability of ``stopping at the top or one unit from the top'' of the random walk (semi best selection problem). This is, the objective is to find the stopping rule $\tau$ that maximizes 
\begin{equation}\label{eq:one}
P(M_N - S_{\tau}\leq 1).
\end{equation}

 In his examination he showed that the strategy is no longer simple, it depends on an interesting sequence of critical values $p_n$ such that when there are $n$ steps remaining until the time horizon, if $M_{N-n} - S_{N-n}=1$ and $p \leq p_n$ it is optimal to stop, otherwise it is optimal to continue. He proved the following properties as well: $p_n<1/2$ for all $n\geq 4$, $\limsup_{n\to\infty} p_n=1/2$ and $p_{2n+1}\geq p_{2n-1}\geq p_{2n}$ for all $n\geq 4$. And further conjectured that (i) $\lim_{n\to\infty} p_n=1/2$ and (ii) $p_{2n}\leq p_{2n+2}$ for all $n\geq 2$. In this work through lengthy calculations the following results are obtained towards proving these conjectures. 

\begin{theorem} We have\\
a. $p_{4}$ is the smallest term of the sequence $p_n$\\  
b. For every $n\geq 0$, $p_{2n+1} \geq p_{2n+4}$  \\
c. For every $n\geq 3$, $p_{2n+1} \geq p_{2n+6}$ 
\end{theorem}

The main purpose of this work is to show the difficulty and the deceiving behavior of $\lim_{n\to\infty} p_n$. See Table $1$ for some values of $p_n$.  
\section{Background for stopping near the top of a random walk}

The content of this section has been extracted from the work by Allaart \cite{Allaart}. The interesting sequence of constants $\{p_{n}\}$is presented in Proposition \ref{prop:PAprop1}. Some of its properties are enclosed in Lemma \ref{eq:lemma2} and Theorem \ref{lem:inequalities}. \\

From now on, the context is to maximize \eqref{eq:one}. Given a finite time horizon $N$, and $0\leq n\leq N$, express
\[M_N - S_n = (M_n-S_n)\vee \max_{n\leq k\leq N}(S_k-S_n)=Z_n \vee M_{N-n}^{'}, \]
where $a\vee b = \max (a,b)$ and $Z_n:=M_n-S_n$. By the independent and stationary increments of the random walk, $M^{'}_{N-n}$ is a random variable independent of the walk up to time $n$, having the same distribution as $M_{N-n}$. Thus, stopping at time $n$ when $Z_n=j$, yields a probability of win $P(j\vee M_k \leq 1)$, where $k=N-n$. This probability is zero if $j\geq 2$, and equals to $P(M_k\leq 1)$ if $j=0$ or $1$.\\

Let the states of the process $(N-n,Z_n)$ be $(k,j)$. The process starts in $(N,0)$ and from $(k,j)$ can move to $(k-1,j-1)$ with probability $p$ or to $(k-1,j+1)$ with probability $q$ if $j>0$. From $(k,0)$ it moves to $(k-1,0)$ with probability $p$ or to $(k-1,1)$ with probability $q$.
 
Note first that in state $(k,j)$ with $j> 1$, and $k>0$ it is trivially optimal to continue. In state $(k,0)$ with $k>0$ it is optimal to continue as well: stopping gives the win probability $P(M_k\leq 1)$, which is at most $P(M_{k-1}\leq 1)$, the win probability of taking one more step and stop. Thus, the critical states are $(k,1)$, $k \in \mathbb{N}$. 
For convenience, we will refer to the critical states defined above as $(n,1)$, $n \in \mathbb{N}$. The following Proposition is about the existence of the interesting sequence of constants.
\begin{proposition}\label{prop:PAprop1}
There exists a sequence $p_n$ of critical values of $p$ such that in state $(n,1)$, it is optimal to stop if and only if $p\leq p_n$.
\end{proposition}
These critical values $p_n$ are zeros of polynomials which can be calculated easily by hand for the first few values of $n$. Table $1$ shows the first $50$ critical values calculated using a code typed in Maple.
\begin{table}[h]
\caption{Values of $p_n$ for $1\leq n \leq 50$, truncated at six decimals.}
\[\begin{array}{c|c|c|c|c|c|c|c|c|c}
\hline
n & p_n & n & p_n & n & p_n & n & p_n &n & p_n \\ \hline
1  & 1           & 11 & 0.484529... & 21 & 0.488177... & 31 & 0.490257... & 41 & 0.491584... \\ 
2  & 0.5         & 12 & 0.479846... & 22 & 0.486398... & 32 & 0.489266... & 42 & 0.490935... \\ 
3  & 0.5         & 13 & 0.485434... & 23 & 0.488683... & 33 & 0.490567... & 43 & 0.491795...  \\ 
4  & 0.468989... & 14 & 0.481753... & 24 & 0.487131... & 34 & 0.489665... & 44 & 0.491191... \\ 
5  & 0.482881... & 15 & 0.486249... & 25 & 0.489137... & 35 & 0.490852... & 45 & 0.491992... \\ 
6  & 0.471448... & 16 & 0.483307... & 26 & 0.487767... & 36 & 0.490027... & 46 & 0.491429... \\ 
7  & 0.482686... & 17 & 0.486970... & 27 & 0.489547... & 37 & 0.491114... & 47 & 0.492176... \\ 
8  & 0.474706... & 18 & 0.484536... & 28 & 0.488327... & 38 & 0.490356... & 48 & 0.491649... \\ 
9  & 0.483572... & 19 & 0.487609... & 29 & 0.489918... & 39 & 0.491358... & 49 & 0.492350... \\ 
10 & 0.477526... & 20 & 0.485545... & 30 & 0.488823... & 40 & 0.490658... & 50 & 0.491855... \\ \hline
\end{array}\]
\end{table}\label{table:t1}

\begin{lemma}\label{eq:lemma2} We have \\
(i) $p_1=1$ and $p_2=p_3=1/2$;\\
(ii) $p_n<1/2$ for all $n\geq 4$; and \\
(iii) $\limsup_{n\to\infty} p_n=1/2$. \\
\end{lemma}

\begin{theorem} \label{lem:inequalities}
$p_{2m+1}\geq p_{2m-1}\geq p_{2m}$ for all $m\geq 4$.
\end{theorem}
A sketch of proof of this Theorem is presented to introduce notation and the necessary formulas for Section $3$. Now, we define the following win probabilities. 
\begin{align*}
V_n:&=\mbox{optimal win probability from state $(n,1)$},\\
W_n:&=\mbox{optimal win probability from state $(n,1)$}\\ &\qquad \mbox{if we continue at least one more step}\\
 &\qquad \mbox{and play optimally from then on},\\
U_n:&=\mbox{win probability if we stop at state $(n,1)$}\\
&=P(M_n\leq 1).
\end{align*}
Then $V_n=\max\{U_n,W_n\}$. Also define $\pi_{n,i}$ to be the optimal win probability from state $(n,i)$ so that in particular, $V_n=\pi_{n,1}$. And a formula for $U_n$ can be derived so that for each $n\in\mathbb{N}$ we obtain 
\begin{equation}
U_{2n}=U_{2n+1}=1-\sum_{j=1}^n t_{j+1}p^{j+1}q^{j-1},
\label{eq:Uformula}
\end{equation}
where $t_m:=\frac{1}{m}\binom{2m-2}{m-1}$. Finally, to sketch the proof consider the formula
%
%
\begin{equation}
P(M_n=k, S_n=l)=a_{n,k,l}p^{(n+l)/2}q^{(n-l)/2}, \qquad 0\leq l\leq k,
\label{eq:formula}
\end{equation}
where
\begin{equation*}
a_{n,k,l}:=\binom{n}{\frac12(n+2k-l)}-\binom{n}{\frac12(n+2k+2-l)}.
\end{equation*}
\textbf{Sketch of Proof of Theorem \ref{lem:inequalities}. } \textit{The proof is reproduced from [1] with the author's permission.} 
Let $m\geq 3$. Let $1/2>p\geq p_{2m+3}$, and assume $p_{2m+3}\geq p_j$ for all $4\leq j\leq 2m+2$. (This may be assumed on account of the induction hypothesis.) If in state $(2m+3,1)$ or $(2m+4,1)$ it is decided to continue, then the optimal strategy is to wait until there are $3$ steps left, and play optimally from then on. 
First it is shown that
\begin{equation}
W_{2m+4}-U_{2m+4}\geq W_{2m+3}-U_{2m+3}.
\label{eq:odd-to-even}
\end{equation}
This inequality implies that, if in state $(2m+3,1)$ it is optimal to continue, then it is optimal to continue in state $(2m+4,1)$ as well; thus, $p_{2m+4}\leq p_{2m+3}$.

Let $\Delta\pi_{3,k}:=\pi_{3,k}-\pi_{3,k+1}$. Then
\begin{equation}
U_{2m+4}-U_{2m+3}=-t_{m+3}p^{m+3}q^{m+1},
\label{eq:U-difference}
\end{equation}
and
\begin{align*}
W_{2m+3}&=\sum_{k=0}^4P(1\vee M_{2m}-S_{2m}=k)\pi_{3,k}\\
&=\sum_{k=0}^4P(1\vee M_{2m}-S_{2m}\leq k)\Delta\pi_{3,k},\\
W_{2m+4}&=\sum_{k=0}^4P(1\vee M_{2m+1}-S_{2m+1}\leq k)\Delta\pi_{3,k},
\end{align*}
so
\begin{equation*}
W_{2m+4}-W_{2m+3}=\sum_{k=0}^4\Delta P_{2m,k}\Delta\pi_{3,k},
\end{equation*}
where
$$\Delta P_{n,k}:=P(1\vee M_{n+1}-S_{n+1}\leq k)-P(1\vee M_n-S_n\leq k).$$
It is easy to see that
\begin{equation*}
P(1\vee M_{n+1}-S_{n+1}\leq k)=pP(M_n-S_n\leq k)+qP(2\vee M_n-S_n\leq k),
\end{equation*}
and since
\begin{gather*}
P(M_n-S_n\leq k)-P(1\vee M_n-S_n\leq k)=P(M_n=0, S_n=-k),\\
P(1\vee M_n-S_n\leq k)-P(2\vee M_n-S_n\leq k)=P(M_n\leq 1, S_n=1-k),
\end{gather*}
it follows that
\begin{equation*}
\Delta P_{n,k}=pP(M_n=0, S_n=-k)-qP(M_n\leq 1, S_n=1-k).
\end{equation*}
Now put $n=2m$ and apply the last identity for $k=0,1,\dots,4$. Also use \eqref{eq:formula} and the notation
$$d_{m,j}:=\binom{2m}{m+j}-\binom{2m}{m+j+1}.$$
Note that $d_{m,j}\geq 0$ for all $m\in\mathbb{N}$ and $j\geq 0$.
Then:
\begin{align*}
\Delta P_{2m,0}&=pP(M_{2m}=0, S_{2m}=0)=d_{m,0}p^{m+1} q^m,\\
\Delta P_{2m,1}&=-qP(M_{2m}\leq 1, S_{2m}=0)=-(d_{m,0}+d_{m,1})p^m q^{m+1},\\
\Delta P_{2m,2}&=pP(M_{2m}=0, S_{2m}=-2)=d_{m,1}p^{m}q^{m+1},\\
\Delta P_{2m,3}&=-qP(M_{2m}\leq 1, S_{2m}=-2)=-(d_{m,1}+d_{m,2})p^{m-1}q^{m+2},\\
\Delta P_{2m,4}&=pP(M_{2m}=0, S_{2m}=-4)=d_{m,2}p^{m-1}q^{m+2}.
\end{align*}

Now (since $p<p_3$), $\pi_{3,0}=1-p^2 q=(p^3+2p^2 q+3p q^2+q^3)$,  $\pi_{3,1}=1-p^2=(2p^2 q+3p q^2+q^3)$, $\pi_{3,2}=(3p^2 q+p q^2)$,  $\pi_{3,3}=p^2=(p^3+p^2 q)$,  $\pi_{3,4}=p^3$, and $\pi_{3,5}=0$. Thus
\begin{align*} 
\Delta\pi_{3,0}&=p^3,\\ 
\Delta\pi_{3,1}&=-p^2 q+2p q^2+q^3,\\
\Delta\pi_{3,2}&=-p^3+2 p^2 q+p q^2,\\
\Delta\pi_{3,3}&=p^2 q,\\
\Delta\pi_{3,4}&=p^3.
\end{align*}

Putting everything together, 
\begin{multline}
W_{2m+4}-W_{2m+3}=p^m q^m[d_{m,0}p^4-d_{m,1}p^3 q+(d_{m,0} +3d_{m,1}+d_{m,2})p^2 q^2\\
\qquad\quad-(2d_{m,0}+2d_{m,1}+d_{m,2})pq^3-(d_{m,0}+d_{m,1})q^4].
\label{eq:W4W3}
\end{multline}
It follows using \eqref{eq:U-difference} that $W_{2m+4}-U_{2m+4}\geq W_{2m+3}-U_{2m+3}$ if and only if
\begin{multline*}
d_{m,0}p^4+(t_{m+3}-d_{m,1})p^3 q+(d_{m,0} +3d_{m,1}+d_{m,2})p^2 q^2\\
-(2d_{m,0}+2d_{m,1}+d_{m,2})pq^3-(d_{m,0}+d_{m,1})q^4\geq 0.
\end{multline*}
This inequality can written in polynomial form and show it is true for $p\leq 1/2$.
 
Next, it is shown that
\begin{equation}
W_{2m+5}-U_{2m+5}\leq W_{2m+3}-U_{2m+3}.
\label{eq:odd-to-odd}
\end{equation}
At $p=p_{2m+3}$, the right hand side of this inequality is zero, so $W_{2m+5}\leq U_{2m+5}$. Thus, \eqref{eq:odd-to-odd} implies that $p_{2m+5}\geq p_{2m+3}$.

Define
$$\Delta^2 P_{n,k}:=P(1\vee M_{n+2}-S_{n+2}\leq k)-P(1\vee M_n-S_n\leq k).$$
Then
\begin{multline*}
P(1\vee M_{n+2}-S_{n+2}\leq k)=p^2P(M_n-S_n\leq k)+2pqP(1\vee M_n-S_n\leq k)\\
+q^2P(3\vee M_n-S_n\leq k).
\end{multline*}
Now
\begin{multline*}
P(1\vee M_n-S_n\leq k)-P(3\vee M_n-S_n\leq k)\\
=P(M_n\leq 1, S_n=1-k)+P(M_n\leq 2, S_n=2-k)
\end{multline*}
and
\begin{equation*}
P(M_n-S_n\leq k)-P(1\vee M_n-S_n\leq k)=P(M_n=0, S_n=-k).
\end{equation*}
Thus,
\begin{multline*}
\Delta^2 P_{n,k}=p^2P(M_n=0, S_n=-k)-q^2P(M_n\leq 1, S_n=1-k)\\
-q^2P(M_n\leq 2, S_n=2-k).
\end{multline*}
Applying this again with $n=2m$ and $k=0,1,\dots,4$, yields
\begin{align*}
\Delta^2 P_{2m,0}&=d_{m,0}p^{m+2}q^m-d_{m,1}p^{m+1}q^{m+1},\\
\Delta^2 P_{2m,1}&=-(d_{m,0}+d_{m,1})p^m q^{m+2},\\
\Delta^2 P_{2m,2}&=d_{m,1}p^{m+1}q^{m+1}-(d_{m,0}+d_{m,1}+d_{m,2})p^m q^{m+2},\\
\Delta^2 P_{2m,3}&=-(d_{m,1}+d_{m,2})p^{m-1}q^{m+3},\\
\Delta^2 P_{2m,4}&=d_{m,2}p^m q^{m+2}-(d_{m,1}+d_{m,2}+d_{m,3})p^{m-1}q^{m+3}.
\end{align*}
Now
\begin{equation*}
U_{2m+5}-U_{2m+3}=-t_{m+3}p^{m+3}q^{m+1}
\end{equation*}
because $U_{2m+5}=U_{2m+4}$. Putting everything together,
\begin{equation}\label{eq:calculation}
\begin{split} 
W_{2m+5}-W_{2m+3}-(U_{2m+5}&-U_{2m+3})\\
&=\sum_{k=0}^4\Delta^2 P_{2m,k}\Delta\pi_{3,k}+t_{m+3}p^{m+3}q^{m+1}\\
&=p^m q^m[(d_{m,0}p^5-2d_{m,1}p^4 q+(d_{m,0}+3d_{m,1}+2d_{m,2})p^3 q^2\\
&-(d_{m,0}+d_{m,1}+3d_{m,2}+d_{m,3})p^2 q^3
-(3d_{m,0}+4d_{m,1}+2d_{m,2})pq^4\\
&-(d_{m,0}+d_{m,1})q^5+ t_{m+3}p^3 q]
\end{split}
\end{equation}
The proof is completed by showing that (\ref{eq:calculation}) is at most zero for $m \geq 3$ and $p$ in the corresponding sub interval of $[0,\frac{1}{2}]$. $\Box$\\
It can be seen in Table $1$ that the even indexes of the sequence are increasing. 
\begin{conjecture}
$p_{2n}\leq p_{2n+2}$ for all $n\geq 2$.
\end{conjecture}
And the main question of this problem.
\begin{conjecture}
$\lim_{n\to\infty} p_n=1/2$.
\end{conjecture}

\begin{section}{ New results}
In this section we present the new results obtained towards proving the conjectures. For further details see \cite{Islas}.\\

\textbf{Proof of Theorem 1.1.a.} It will be proven by induction. By direct calculation of $p_{1}$, $p_{2}$, ..., $p_{27}$, the statement holds for $n=1,2,...,27$. Suppose $p_k\geq p_4$ for $1\leq k\leq n-1$. It must be shown that  for $p \leq p_{4}$
\begin{equation}
  U_{n} - W_{n} \geq 0 
\end{equation}\label{eq:ineqprop3}
Let $n\geq 2$, and define the stopping time
\begin{equation}\label{eq:sigma}
\sigma:=\inf\{j\geq 1: 1\vee M_j-S_j=1\},
\end{equation}
which is the first time the random walk from state $(n,1)$ comes back one unit below its running maximum. By the induction hypothesis, the optimal strategy from state $(n,1)$ is to stop at each critical state. Thus the win probability is  
\[
W_n = \sum_{j=2}^n P(\sigma=j)U_{n-j} + p^n,
\]
since if in state $(n,1)$ it is decided to continue, there can be a win only if the walk either records a string of $n$ straight up-steps or comes back to one unit below its running maximum at some future time. And it was shown in the proof of \ref{prop:PAprop1} (ommited) that
\[
U_n = \sum_{j=2}^n P(\sigma=j)U_{n-j} - \sum_{j=3}^n p^{j-1}P(M_{n-j+1}=0) +qP(\tau_1>n-1).  
\]
where the hitting times $\tau_j:=\inf\{n\geq 0:S_n=j\}$. It follows that  
\begin{equation}
U_n - W_n = - \sum_{j=3}^n p^{j-1}P(M_{n-j+1}=0) +qP(\tau_1>n-1) - p^n\\
=a_n - \sum_{i=1}^{n-1} a_{n-i}P(\tau_{1}=i)
\end{equation}
where $a_j =1 - \frac{p}{q}(1-p^j)$ and it is decreasing in $j$.\\
First, an upper bound for
\begin{equation}
 \sum_{j=1}^{n-1} a_{n-i}P(\tau_{1}=i)
\end{equation}\label{eq:sumaprop3}
is found using the formula \begin{equation*}
P(\tau_1=2j+1)=\frac{1}{j+1}\binom{2j}{j}p^{j+1}q^j.
\end{equation*}
Suppose $n$ is even, then
\begin{align*}
\sum_{i=1}^{n-1} a_{n-i}P(\tau_{1}=i) &= \sum_{i=1}^{n-7} a_{n-i}P(\tau_{1}=i) + a_5P(\tau_1=n-5)+ a_3P(\tau_1=n-3)+ a_1P(\tau_1=n-1) \\
&\leq a_7\sum_{i=1}^{\infty} P(\tau_{1}=i) + a_5P(\tau_1=n-5)+ a_3P(\tau_1=n-3)+ a_1P(\tau_1=n-1).
\end{align*}
The summation equals $p/q$, thus
\begin{align*}
\sum_{i=1}^{n-1} a_{n-i}P(\tau_{1}=i) &\leq a_7 \frac{p}{q} + a_5\frac{1}{j_5+1} \binom{2j_5}{j_5}p^{j_5+1}q^{j_5}+ a_3\frac{1}{j_3+1} \binom{2j_3}{j_3}p^{j_3+1}q^{j_3}\\
&+a_1\frac{1}{j_1+1} \binom{2j_1}{j_1}p^{j_1+1}q^{j_1}
\end{align*}
where 
\[j_k:=\frac{n-k-1}{2} , \qquad k \in \mathbb{N}.\] 

The binomial coefficients are bounded using Stirling's approximation for factorials 
\[
\left( 2\pi\right)^{1/2}j^{j+1/2}e^{-j} \leq j!\le e j^{j+1/2}e^{-j}
\]
and the inequality $e \geq \left( 2\pi\right)^{1/2}e^{1/(12j)}$ to obtain
\begin{equation*}
\frac{1}{j+1} \binom{2j}{j} \leq \frac{4^je}{(2j)^{1/2}\pi(j+1)} \leq \frac{4^je}{2^{1/2}\pi j^{3/2}}.
\end{equation*}
From the above inequality and $p^{j+1}q^{j}\le p4^{-j}$ it follows that 
\begin{align*}
\sum_{j=1}^{n-1} a_{n-i}P(\tau_{1}=i) &\leq a_7 \frac{p}{q} a_5 \frac{4^{j_5}e}{2^{1/2}\pi j_5^{3/2}}p4^{-j_5} +  a_3 \frac{4^{j_3}e}{2^{1/2}\pi j_3^{3/2}}p4^{-j_3}+ a_1 \frac{4^{j_1}e}{2^{1/2}\pi j_1^{3/2}}p4^{-j_1}\\
& \leq a_7 \frac{p}{q} + \frac{ep}{2^{1/2}\pi j_5^{3/2}}(a_5+a_3+a_1). 
\end{align*}
Similarly, when $n$ is odd 
\begin{align*}
\sum_{i=1}^{n-1} a_{n-i}P(\tau_{1}=i) &\le  a_8\sum_{i=1}^{n-8} P(\tau_{1}=i) + a_6P(\tau_1=n-6)+ a_4P(\tau_1=n-4)+ a_2P(\tau_1=n-2)\\
&\le a_8 \frac{p}{q} + \frac{ep}{2^{1/2}\pi j_6^{3/2}}(a_6+a_4+a_2). 
\end{align*}
Then for all $n$, recalling that $a_j$ is decreasing in $j$,
\begin{align*}
\sum_{i=1}^{n-1} a_{n-i}P(\tau_{1}=i) &\le a_7 \frac{p}{q} + \frac{ep}{2^{1/2}\pi j_6^{3/2}}(a_5+a_3+a_1)\\
&= a_7 \frac{p}{q} + \frac{ep}{2^{1/2}\pi (\frac{n-7}{2})^{3/2}}(a_5+a_3+a_1). 
\end{align*}
Since the expression on the right decreases in $n$, it will be enough to use $n=27$. For $n \geq 27$,
\[
\sum_{i=1}^{n-1} a_{n-i}P(\tau_{1}=i) \le a_7 \frac{p}{q} + \frac{p}{50}(a_5+a_3+a_1).  
\]
The desired bound has been found, thus 
\begin{align*}
U_n-W_n&= a_n - \sum_{i=1}^{n-1} a_{n-i}P(\tau_{1}=i)\\
&\ge a_n - a_7 \frac{p}{q} - \frac{p}{50}(a_5+a_3+a_1) \\
&\geq 1 - \frac{p}{q} - a_7 \frac{p}{q} - \frac{p}{50}(a_5+a_3+a_1) \\
&=\frac{1}{50}\frac{(50p^8+49p^7+50p^6+49p^5+50p^4+49p^3+56p^2-153p+50)}{q}. (*)
\end{align*}
It now remains to show that the polynomial in the numerator is positive on $[0,p_4]$. Its derivative
$400p^7+343p^6+300p^5+245p^4+200p^3+147p^2+112p-153$ is clearly an increasing function on $[0,1/2]$ and at the end point $1/2$ the function takes the value -$2.0781250$. Hence the numerator in (*) is a decreasing function. For $p=.47>p_4$ the numerator is $0.01762055$.\\
\textbf{Proof of Theorem 1.1.b.} It must be proven that $p_{2m+1}\geq p_{2m+4}$ for all $m\geq 0$. 
Note that by direct calculation of $p_1,p_2,...,p_8$ the statement is true for $m= 0,1,2$. Let $m> 2$. It will be sufficient to show that for $p_{2m+3}> p\geq p_{2m+1}$, 
  \begin{equation} \label{eq:ineqprop1}
  W_{2m+1} - U_{2m+1} \leq W_{2m+4} - U_{2m+4} 
  \end{equation}
holds. The previous inequality implies that if in state $(2m+1,1)$ it is optimal to continue, then it's optimal to continue in state $(2m+4,1)$ as well; thus $p_{2m+1} \geq p_{2m+4}$. By Theorem \ref{lem:inequalities},  $p_{2m+1}\geq p_j$ for all $4\leq j\leq 2m+2$. Then, the optimal strategy from state $(2m+1,1)$ is to wait until $3$ steps are left and then play optimally from then on. From state $(2m+4,1)$, the first critical state one can reach is $(2m+2,1)$, in which it is optimal to continue $(p\geq p_{2m+1}\geq p_{2m+2})$, hence the optimal strategy is to wait until $3$ steps are left and then play optimally from then on. 
To show \eqref{eq:ineqprop1}, add 
\[-W_{2m+3} + U_{2m+3} \]
to both sides, so  \eqref{eq:ineqprop1} holds if and only if 
\begin{equation} \label{eq:ineq2prop1}
   \underbrace{(W_{2m+4} - W_{2m+3} - U_{2m+4} + U_{2m+3})}_{L} + \underbrace{W_{2m+3} -W_{2m+1} - (U_{2m+3} - U_{2m+1})}_{R} \geq 0
\end{equation}
L is found using \eqref{eq:U-difference} and \eqref{eq:W4W3}:
\begin{align*} 
L&= p^m q^m[d_{m,0}p^4-d_{m,1}p^3 q+(d_{m,0} +3d_{m,1}+d_{m,2})p^2 q^2\\
&-(2d_{m,0}+2d_{m,1}+d_{m,2})pq^3-(d_{m,0}+d_{m,1})q^4] + t_{m+3}p^{m+3}q^{m+1},
\end{align*}
and R is found replacing $m$ by $(m-1)$ in \eqref{eq:calculation} 
\begin{align*} 
R&= p^{m-1} q^{m-1} [ d_{m-1,0}p^5 - 2d_{m-1,1}p^4q + (d_{m-1,0}+3d_{m-1,1}+ 2d_{m-1,2})p^3q^2\\
&- (d_{m-1,0}+d_{m-1,1}+3d_{m-1,2}+d_{m-1,3})p^2q^3 - (3d_{m-1,0}+4d_{m-1,1}+2d_{m-1,2})pq^4 \\
&- (d_{m-1,0}+d_{m-1,1})q^5 + t_{m+2}p^3q ].
\end{align*}
After long and tedious calculations it follows that   
\begin{equation}
L+R\geq 0
\end{equation}\label{eq:LR}
if and only if
\[
pol:=m^{6}b_{6} + m^{5}b_{5} + m^{4}b_{4} + m^{3}b_{3} + m^{2}b_{2} + m b_{1}\geq 0,
\]
where
\begin{align*}
b_6 =& -48p^6 + 240p^5 - 340p^4 + 112p^3 + 62p^2 -21p -4, \\
b_5 =& -192p^6 + 864p^5 - 1014p^4 + 468p^2 -87p -30, \\
b_4 =& -228p^6 + 852p^5 -652p^4 - 572p^3 + 722p^2 -15p -76, \\
b_3 =&  -48p^6 + 72p^5 + 78p^4 -72p^3 -108p^2 + 195p -66,\\
b_2 =&  60p^6 -228p^5 + 128p^4 + 316p^3 -280p^2 + 36p +8,\\
b_1 =&  24p^6 -72p^5 + 72p^4 -72p^3+ 144p^2 - 108p + 24.\\
\end{align*}
Minimizing the previous polynomials numerically on the interval $0.478\leq p\leq 0.5$ gives 
\begin{align*}
  \min(b_6) &\approx  0.02685526, \\
	\min(b_5) &\approx 1.67875176, \\
  \min(b_4) &\approx 3.82783882, \\
  \min(b_3) &\approx -0.03365067, \\
  \min(b_2) &\approx -2.687500, \\
  \min(b_1) &\approx -0.375. 
\end{align*}
Thus 
\begin{align*}
pol\geq 0.0268m^{6}+1.6787m^{5} + 3.8278m^{4} -0.0336m^3 -2.6875m^2 -0.375m.
\end{align*}
The last polynomial has six real zeros; the unique positive zero is $m=0.7860814063$. At $m=1$ it is positive, therefore it remains positive for all $m\geq 1$. It follows that \eqref{eq:ineqprop1} is true for $m> 2$ and $p_{2m+3}> p \geq p_{2m+1}$, since $p_{2m+1}\geq 0.482881151>0.478$ for all $m\geq 2$.\\   
Now, consider the following win probability. 
\begin{center}
$\widetilde{W}_{n}:=$ win probability from state $(n, 1)$ using the stopping rule $\rho$,\\
\end{center}
where $\rho$ is the stopping rule that waits until there are $3$ steps left and plays optimally from then on. \\
\textbf{Proof of Theorem 1.1.c.} Note that by direct calculation of $p_1,p_2,...,p_{46}$ the statement is true for all $m\leq20$. Let $m> 20$. It will be sufficient to show that for $p_{2m+3}> p \geq p_{2m+1}$ 
  \begin{equation} \label{eq:ineqprop2}
  W_{2m+1} - U_{2m+1} \leq W_{2m+6} - U_{2m+6} 
  \end{equation}
holds. The previous inequality implies that if in state $(2m+1,1)$ is optimal to continue, then it's optimal to continue in state $(2m+6,1)$ as well; thus $p_{2m+1} \geq p_{2m+6}$. By Theorem \ref{lem:inequalities}, $p_{2m+1} \geq p_j$ for all $4\leq j \leq 2m+2$. Then, the optimal strategy from state $(2m+1,1)$ is to wait until $3$ steps are left and then play optimally from then on. From state $(2m+6,1)$, it will be enough to use the strategy that waits until $3$ steps are left and plays optimally from then on; the win probability with this strategy is $\widetilde{W}_{2m+6}$. 
To show  \eqref{eq:ineqprop2}, add
\[\widetilde{W}_{2m+5} + W_{2m+3}\]
 to both sides, so  \eqref{eq:ineqprop2} holds if and only if
\begin{equation*} 
  ( W_{2m+6} - \widetilde{W}_{2m+5}) + (U_{2m+1} - U_{2m+6}) + (\widetilde{W}_{2m+5} - W_{2m+3}) + (W_{2m+3}-W_{2m+1})  \geq 0.
\end{equation*}
Using the fact that $W_{n} \geq \widetilde{W}_{n}$, it is enough to show that 
\begin{equation} \label{eq:ineq2prop2}
   \underbrace{\widetilde{W}_{2m+6} - \widetilde{W}_{2m+5}}_{a} + \underbrace{U_{2m+1} - U_{2m+6}}_{b}  + \underbrace{\widetilde{W}_{2m+5} - W_{2m+3}}_{c}  + \underbrace{W_{2m+3}-W_{2m+1}}_{d} \geq 0.
\end{equation}
To obtain $\widetilde{W}_{2m+6} - \widetilde{W}_{2m+5}$ replace $(m)$ by $(m+1)$ in \eqref{eq:W4W3}, calculate $U_{2m+1} - U_{2m+6}$ by direct use of \eqref{eq:Uformula}, for $\widetilde{W}_{2m+5} - W_{2m+3}$   
subtract $t_{m+3}p^{m+3}q^m$ from \eqref{eq:calculation}, and use $m-1$ instead of $m$ to get $W_{2m+3}-W_{2m+1}$, then
\begin{align*}
\widetilde{W}_{2m+6} - \widetilde{W}_{2m+5}&= p^{m+1} q^{m+1}[d_{m+1,0}p^4-d_{m+1,1}p^3 q+(d_{m+1,0} +3d_{m+1,1}+d_{m+1,2})p^2 q^2\\
&-(2d_{m+1,0}+2d_{m+1,1}+d_{m+1,2})pq^3-(d_{m+1,0}+d_{m+1,1})q^4], \\
U_{2m+1} - U_{2m+6} &= t_{m+2}p^{m+2}q^m +  t_{m+3}p^{m+3}q^{m+1} + t_{m+4}p^{m+4}q^{m+2}, \\
\widetilde{W}_{2m+5}-W_{2m+3} &=p^m q^m [(d_{m,0}p^5-2d_{m,1}p^4 q+(d_{m,0}+3d_{m,1}+2d_{m,2})p^3 q^2\\
&\qquad-(d_{m,0}+d_{m,1}+3d_{m,2}+d_{m,3})p^2 q^3
-(3d_{m,0}+4d_{m,1}+2d_{m,2})pq^4\\
&-(d_{m,0}+d_{m,1})q^5 ]\\
W_{2m+3}-W_{2m+1} &=p^{m-1} q^{m-1} [(d_{m-1,0}p^5-2d_{m-1,1}p^4 q+(d_{m-1,0}+3d_{m-1,1}+2d_{m-1,2})p^3 q^2\\
&-(d_{m-1,0}+d_{m-1,1}+3d_{m-1,2}+d_{m-1,3})p^2 q^3
-(3d_{m-1,0}+4d_{m-1,1}+2d_{m-1,2})pq^4\\
&-(d_{m-1,0}+d_{m-1,1})q^5 ].
\end{align*}
Long and tedious calculations show that \eqref{eq:ineq2prop2} is true  if and only if\\
\[
pol:= m^{6}b_{6} + m^{5}b_{5} + m^{4}b_{4} + m^{3}b_{3} + m^{2}b_{2} + m b_{1} + b_{0}\geq 0, 
\]
where
\begin{align*}
b_6 =& 192p^8-1152p^7+2320p^6-1776p^5+140p^4+336p^3-34p^2-21p-4,\\
b_5 =& 1248p^8-7104p^7+13312p^6-8352p^5-1702p^4+2880p^3-52p^2-171p-46, \\
b_4 =& 2928p^8-15840p^7+27724p^6-14100p^5-7324p^4+6636p^3+602p^2-363p-196, \\
b_3 =& 2760p^8-14160p^7+23200p^6-9720p^5-7690p^4+5160p^3+860p^2+135p-370, \\
b_2 =& 408p^8-1728p^7+2140p^6-204p^5-1240p^4+1284p^3-976p^2+816p-256, \\
b_1 =& -768p^8+3984p^7-6632p^6+2952p^5+1832p^4-840p^3-448p^2+36p+56, \\
b_0 =& -288p^8+1440p^7-2304p^6+960p^5+864p^4-1056p^3+768p^2-432p+96. \\
\end{align*}
Minimizing the polynomials with respect to $p$ over $[0.4925,0.5]$ gives
\begin{align*}
  \min(b_6) &\approx  0.000507117,\\
  \min(b_5) &\approx 3.43282349, \\
  \min(b_4) &\approx 18.9614840, \\
  \min(b_3) &\approx 28.5306405, \\
  \min(b_2) &\approx 3.4720812,\\
  \min(b_1) &= -11.75, \\
  \min(b_0) &= -1.875. 
\end{align*}
Thus,
\begin{align*}
pol\geq 0.0005m^6 + 3.4328m^5 + 18.9614m^4 + 3.4720m^3  -11.75m^2  + m  -1.875
\end{align*}
and this lower bound is clearly positive for $m\geq 1$. Thus, for each $m> 20$ and $p_{2m+3}> p\geq p_{2m+1}$ \eqref{eq:ineqprop2} is true. $\Box$ 

\section{Conclusions}
Further computations show that the pattern $p_{2m+8 }\geq p_{2m+1}$ holds for $5\leq m \leq 156$ but for a greater $m$, the inequality is reversed making the sequence tricky and more difficult to work with. The intriguing oscillation pattern of the sequence, the unsolved conjecture for the even critical values, and the computational difficulty increased with $n$, make useless the effort to try to solve the main conjecture. 
\section*{Acknowledgment}
To Pieter Allaart for his permission to reproduce the proof of Theorem $2.3$.
\end{section}


\footnotesize

\begin{thebibliography}{5}
	
\bibitem{Allaart}
{\sc  P. C. Allaart}, How to stop near the top in a random walk? In {\em Decision Making Processes under Uncertainty and Ambiguity} RIMS Kokyuroku, pp. 33-40, 2010.

\bibitem{Hlynka}
{\sc M. Hlynka} and {\sc J.N. Sheahan}, The secretary problem for a random walk. In {\em Stoch. Proc. Appl.} 28, pp. 317-325, 1998.

\bibitem{Yam}
{\sc S. C. P. Yam, S. P. Yung}, and {\sc W. Zhou}, Two rationales behind `buy-and-hold or sell-at-once'. In {\em  J. Appl. Probab.} 46, pp. 651-668, 2009.

\bibitem{Islas}
{\sc J. A. Islas Anguiano}, Optimal Strategies for Stopping Near the Top of a Sequence, dissertation, December 2015; Denton, Texas. (digital.library.unt.edu/ark:/67531/metadc822812/: accessed November 20, 2017), University of North Texas Libraries, Digital Library, digital.library.unt.edu. 



\end{thebibliography}
\end{document}